\documentclass[12pt,twoside,leqno]{article}

\usepackage{amssymb}
\nonstopmode

\def\epv{{$\mbox{}$\hfill ${\Box}$\vspace*{1.5ex} }}
\def\vsp{\vspace*{1.5ex}}
\def\prin{\mbox{\rm prin}}

\def\mod{\mbox{\rm mod}}

\def\modsp{{  \mbox{\rm mod} _ {\mbox{\scriptsize sp} }}}
\def\Ext{\mbox{\rm Ext}}
\def\Hom{\mbox{\rm Hom}}
\def\Ker{\mbox{\rm Ker}}
\def\Im{\mbox{\rm Im}}

\def\Aut{\mbox{\rm Aut}}
\def\id{\mbox{\rm id}}
\def\soc{\mbox{\rm soc}}
\def\bdim{\mbox{\bf dim}}
\def\longarr#1#2{{\buildrel{#1} \over {\hbox to #2pt{\rightarrowfill}}}}
\def\mapdown#1{\Big\downarrow\rlap{$\vcenter{\hbox{$\scriptstyle#1$}}$}}

\newcommand{\CB}{{\mathcal{B}}}

\newcommand{\CH}{{\cal H}}

\newcommand{\CS}{{\cal S}}

\def\ya{\mbox{\em (a) \hspace{0.3em}}}
\def\yb{\mbox{\em (b) \hspace{0.3em}}}
\def\yc{\mbox{\em (c) \hspace{0.3em}}}

\begin{document}

\renewcommand{\thefootnote}{}
\def\refname{\begin{center}\normalsize{\it REFERENCES}\end{center}}

 \title{The existence of Hall polynomials for posets of finite prinjective type}

\author{ Justyna Kosakowska\thanks{Partially supported by Polish KBN Grant 5 P03A 015 21} \\  \it Faculty of Mathematics and Computer Science,\\
\it Nicolaus Copernicus University,\\ \it ul. Chopina 12/18,
87-100 Toru$\acute{n}$, Poland,
\\ \it e-mail justus@mat.uni.torun.pl} \date{}\maketitle



\hfill Dedicated to Professor Daniel Simson

\hfill on the occasion of his 65th birthday \vsp \vsp\\



\begin{abstract} We prove the~existence of Hall polynomials for prinjective
representations of finite partially ordered sets of finite
prinjective type. In Section 4 we shortly discuss consequences of
the existence of Hall polynomials, in particular, we are able to
define a~generic Ringel-Hall algebra for prinjective
representations of posets of finite prinjective type.
\end{abstract}

\section{Introduction} Let $K$ be a~finite field and let $A$ be
a~finite dimensional associative, basic $K$-algebra. All modules
considered in the present paper are right, finite dimensional
$A$-modules. Given $A$-modules $X,Y,Z$, denote by $F^Y_{Z,X}$ the
number of submodules $U\subseteq Y$ such  that $U\simeq X$ and
$Y/U\simeq Z$. Moreover denote by $\Gamma_A$ the Auslander-Reiten
quiver of the algebra $A$. The reader is referred to \cite{ars},
\cite{ass} and to \cite{ri} for the definitions and the
introduction to the theory of representations of algebras.

Let $\Gamma=(\Gamma_0,\Gamma_1)$ be a~directed Auslander-Reiten quiver, with
the set of vertices $\Gamma_0$ and set of arrows $\Gamma_1$. Recall that
for any field $K$ and any $K$-algebra $A$ such that
$\Gamma_A=\Gamma$, we may identify a~function $a:\Gamma_0\to
\mathbb{N}$ with the corresponding $A$-module $M(A,a)=M(a)$ (see
\cite{ringel1}). It was proved by C. M. Ringel (in \cite{ringel1})
that for any directed  Auslander-Reiten quiver $\Gamma$ and all
functions $a,b,c:\Gamma_0\to \mathbb{N}$, there exist polynomials $\varphi^b_{ca}\in
\mathbb{Z}[T]$ with the following property: if $K$ is a~finite
field, and $A$ a~$K$-algebra with $\Gamma_A=\Gamma$ and
symmetrization index $r$, then
$F^{M(A,b)}_{M(A,c),M(A,a)}=\varphi^b_{ca}(|K|^r)$. The
polynomials $\varphi^b_{ca}$ are called {\bf Hall polynomials}.
Moreover, in \cite{ringel92}, C. M. Ringel conjectured the
existence of Hall polynomials for every representation finite
algebra. In \cite{peng} it was proved that there exist Hall
polynomials for representation-finite trivial extension algebras.
The existence of Hall polynomials for cyclic symmetric algebras
was proved in \cite{guo}.

Now we present consequences of the existence of Hall polynomials.
We restrict our considerations to hereditary algebras. Let
$\Delta$ be a~Dynkin quiver, $A=K\Delta$ -- path algebra of
$\Delta$ and $q\in \mathbb{C}$. Following \cite{ringel} we define
$\CH_q(\Delta)$ to be the free abelian group with basis
$(u_M)_{[M]}$, indexed by the set of isomorphism classes of finite
dimensional right $A$-modules. $\CH_q(\Delta)$ is an~associative
ring with identity $u_0$, where the multiplication is defined by
the formula
$$u_{X_1}u_{X_2}=\sum_{[X]}\varphi^X_{X_1,X_2}(q)u_X,$$ and sum runs over
all isomorphism classes of $A$-modules. We call $\CH_q(\Delta)$
{\bf the Ringel-Hall algebra} of $A$.

The motivation for the study of Hall polynomials and Hall algebras
comes from their connection with generic extensions, Lie algebras
and quantum groups (see \cite{ringel}, \cite{ringel1},
\cite{ringel92}, \cite{reineke2001}). It is known that
$\CH_1(\Delta)\otimes_{\mathbb{Z}}\mathbb{C}$ is isomorphic with
the universal enveloping algebra $U({\bf n}_+)$ of ${\bf n}_+$,
where ${\bf g}={\bf n}_-\oplus {\bf h}\oplus {\bf n}_+$ is
a~triangular decomposition of the semisimple complex Lie algebra
${\bf g}$ of type $\Delta$ (see \cite{ringel}).

In the present paper we are interested in an~analogous problem of
the~existence of Hall polynomials for prinjective modules over
incidence algebras of posets of finite prinjective type (see
Section 2 for definition). We define also (Section 4) prinjective
Ringel-Hall algebras for such posets. The paper is organised as
follows. In Section 2 we prove some results concerning injective
and surjective homomorphisms between prinjective modules and we
recall main definitions and results concerning prinjective
modules. In Section 3 the existence of Hall polynomials for
prinjective representations of posets of finite prinjective type
is proved. Section 4 contains consequences of the existence of
Hall polynomials. In particular we give there a~definition of
prinjective Ringel-Hall algebra. Concluding remarks are also
presented in Section 4.

The motivations for the study of prinjective $KI$-modules is the
fact that many of the representation theory problems can be
reduced to the corresponding problems for poset representations
and prinjective modules (see \cite{arn}, \cite{ri}, \cite{s92},
\cite{s93}, \cite{s95}). Prinjective $KI$-modules play
an~important role in the representation theory of finite
dimensional algebras (see \cite{ri}, \cite[Chapter 17]{s92}) and
lattices over orders (see \cite[Chapter 13]{s92}, \cite{s93},
\cite{s95}, \cite{s97a}). Moreover the study of prinjective
modules is equivalent to the study of a~class of bimodule matrix
problems in the sense of Drozd (see \cite{ps}, \cite[Chapter
17]{s92}). \vsp

{\bf Acknowledgment.} The author would like to thank S. Kasjan for
careful reading of this paper and helpful remarks. The main
results of this paper were presented on the X ICRA in Patzcuaro
(Mexico) 2004, on the seminar in Bielefeld (Germany) 2004 and on
the NWDR Workshop in Muenster (3th December 2004) during the stay
of the author supported by Lie Grits (C 0105704).

\vsp

\section{Counting surjective homomorphisms} Let $I=(I,\preceq)$ be
a~finite poset (i.e. partially ordered set) with the partial order
$\preceq$. Let $\max\, I$ denote the set of all maximal elements
of $I$ and $I^-=I\setminus \max\, I$. Given a~field $K$ we denote
by $KI$ the {\bf incidence} $K${\bf -algebra} of the poset $I$,
that is, $$KI=\{(\lambda_{ij})\in \mathbb{M}_I(K)\; ;\;\;
\lambda_{ij}=0 \mbox{ if } i\npreceq j\mbox{ in } I\} \subseteq
\mathbb{M} _I(K)$$  (see \cite{s92}, \cite{s93}). The reader is
referred to \cite{s92}, \cite{s93}, \cite{s95}, \cite{s97a} for
a~discussion of incidence algebras and their applications to the
integral representation theory. A~$KI$-module $X$ may be
identified with the representation $(X_i,\varphi_{ij})_{i\preceq
j\in I}$ of the poset $I$ (i.e. $X_i$ is a~$K$-vector space for
any $i\in I$ and, for all relations $i\preceq j$ in $I$,
$\varphi_{ij};X_i\to X_j$ are linear maps satisfying
$\varphi_{jk}\varphi_{ij}=\varphi_{ik}$ if $i\preceq j\preceq k$).
Recall that
 the dimension vector $\bdim\, X\in \mathbb{Z}^I$ of $X$ is defined by
 $(\bdim\,X)(i)=\dim_K X_i$ for all $i\in I$.
Denote by $P(i)$ the projective $KI$-module corresponding to the
vertex $i$. Without loss of generality we may
assume that $I\subseteq \mathbb{N}$ and that the order $\preceq$
in $I$ is such that $i\preceq j$ in $I$ implies $i\leq j$ in the
natural order. In this case the algebra $KI$ has the following {\bf bipartition}
$$KI=\left[\begin{array}{cc}KI^-&M\\ 0&K(\max\,
I)\end{array}\right], \leqno(2.1)$$ where $M$ is
a~$KI^-$-$K(\max\,I)$-bimodule.

It is well-known (see \cite{s90},\cite[III.2]{ars}) that
a~finitely generated $KI$-module $X$ may be also identified with
the triple
$$X=(X',X'',\varphi:X'\otimes_{KI^-}M\to X''),$$
where $X'$ is a~$KI^-$-module, $X''$ is a~$K(\max\, I)$-module and $\varphi$ is
a~$K(\max\,I)$-module homomorphism. A~homomorphism $f:X\to Y=(Y',Y'',\psi)$
of $KI$-modules is identified with a~pair $(f',f'')$, where $f':X'\to Y'$
is a~$KI^-$-module homomorphism, $f'':X''\to Y''$ is a~$K(\max\, I)$-module homomorphism and
$f''\varphi=\psi (f'\otimes \id)$.
Equivalently, we may identify $X$ with the triple
$$X=(X',X'',\overline{\varphi}:X'\to\Hom_{K(\max\, I)}(M,X'')),$$
where $X'$ is a~$KI$-module, $X''$ is a~$KI^-$-module and
$\overline{\varphi}$ is the~$KI^-$-module homomorphism adjoint to
$\varphi$. A~homomorphism $f:X\to Y=(Y',Y'',\psi)$ of
$KI$-modules, in this case, is identified with a~pair $(f',f'')$,
where $f':X'\to Y'$ is a~$KI^-$-module homomorphism, $f'':X''\to
Y''$ is a~$K(\max\, I)$-module homomorphism and $\overline{\psi}
f'=\Hom_B(M,f'')\overline{\varphi}$. In the present paper we use
and need these three presentations of a~$KI$-module $X$.

Let $\mod(KI)$
denotes the category of all finite dimensional right $KI$-modules.

A~$KI$-module $X$ is said to be {\bf prinjective} if the
$KI^-$-module $X'$ is projective. Let us denote by $\prin(KI)$ the
full subcategory of $\mod(KI)$ consisting of prinjective
$KI$-modules. Note that any projective $KI$-module is prinjective.
The algebra $KI$ is said to be of {\bf finite prinjective type} if
the category $\prin(KI)$ is of finite representation type, i.e.
there exist only finitely many isomorphism classes of
indecomposable prinjective $KI$-modules. \vsp

{\sc Remark.} If the poset $I$ is of finite prinjective type, the
$K$-algebra $KI$ may be of infinite representation type (even
wild). Moreover the category of prinjective modules is not closed
under submodules. Therefore the problem of the~existence of Hall
polynomials for prinjective modules does not reduce to the
corresponding one for representation directed algebras and
Ringel's arguments given in \cite{ringel1} does not apply directly
in our case. In this section wee present a~reduction which allows
us, in Section 3, to develop Ringel's arguments in our case.\epv

Let us denote by $\modsp(KI)$ the full subcategory of $\mod(KI)$
consisting of socle projective modules, i.e. modules $X$ which
have projective socle $\soc(X)$. Following \cite{s90} we define
the functor
$$\Theta:\prin(KI)\to\modsp(KI)$$
by $$(X',X'',\varphi)\mapsto (\Im \overline{\varphi},X'',
\overline{j_\varphi})=(\Theta(X'),\Theta(X)'',\overline{j_\varphi}),$$
where $\overline{j_\varphi}$ is the adjoint map to the inclusion
$j_\varphi:\Im \overline{\varphi}\hookrightarrow \Hom_{K(\max\, I
)}(M,X'')$. Let us collect some properties of these categories and
functor. \vsp

{\sc Lemma 2.2.} {\it  \ya A~$KI$-module $X=(X',X'',\varphi)$ belongs to the category $\modsp(KI)$ if and only if
$\soc(X)$ has the form $(0,Y,0)$, where $Y$ is a~$K(\max\, I)$-module.

\yb The functor $\Theta$ is full and dense with $\Ker\, \Theta
=[(P,0,0)\, ;\; P\mbox{ {\rm projective} } \allowbreak
KI^-\mbox{\rm -module}]$. Moreover $\Theta$ establishes
a~bijection between indecomposable modules which are not in
$\Ker\,\Theta$ and indecomposable modules in $\modsp(KI)$. } \vsp

{\bf Proof.} See \cite{s90} and \cite{ps}.\epv

Now we prove some facts about surjective and injective
homomorphisms of $KI$-modules. These facts are essentially used in
Section 3.\vsp

{\sc Lemma 2.3.} {\it \ya Let $X=(X',X'',\overline{\varphi})$,
$Y=(Y',Y'',\overline{\psi})$ be modules in $\prin(KI)$  and let
$f=(f',f''):X\to Y$ be an~injective (resp. surjective)
$KI$-homomorphism. Then $\Theta(f)$ is an~injective (resp.
surjective) $KI$-homomorphism.

\yb Let $X=(X',X'',\overline{\varphi})$,
$Y=(Y',Y'',\overline{\psi})$ be modules in $\prin(KI)$ and let
$f:X\to Y$ be a~$KI$-homomorphism such that
$\Theta(f)=(g',g''):\Theta(X)\to \Theta(Y)$ is surjective. If $Y$
has no direct summand of the form $(P,0,0)$, where $P$ is
a~projective $KI^-$-module, then $f$ is surjective. }\vsp

{\bf Proof.} (a) Let $f:X\to Y$ be a~homomorphism and
$$g=(g',g'')=\Theta(f)=(\Hom_{K(\max\,I)}(M,f'')|_{\Im\overline{\varphi}},f'').$$
Assume that $f$ is injective. Then the morphisms $f'$ and
$f''=g''$ are injective. Note that $g'$ is injective, because
$f''$ is injective and the functor $\Hom_{K(\max\, I)}(M,-)$ is
left exact.

Now let $f$ be surjective. Then $f'$, $f''=g''$ are surjective. We
have to show that $g':\Im \overline{\varphi}\to \Im
\overline{\psi} $ is surjective. Note that $$g'(\Im
\overline{\varphi})= \Hom_{K(\max\, I
)}(M,f'')(\Im\overline{\varphi})=\Hom_{K(\max\, I
)}(M,f'')\overline{\varphi}(X')=\overline{\psi} f'(X').$$ Since
$f'$ is surjective we have $$\overline{\psi}
f'(X')=\overline{\psi} (Y')=\Im\overline{\psi}.$$ Therefore $g'$
and $g$ are surjective. This finishes the proof of (a).

(b) Let $X=(X',X'',\overline{\varphi})$,
$Y=(Y',Y'',\overline{\psi})$ be modules in $\prin(KI)$  and let
$\Theta(f)=(\Hom_{K(\max\,
I)}(M,f'')|_{\Im\overline{\varphi}},f'')=(g',g''):\Theta(X)\to
\Theta(Y)$ be surjective. It follows that
 $g'\overline{\varphi}:X'\to \Theta(Y)'$
and $\overline{\psi}f'=g'\overline{\varphi}:X'\to \Theta(Y)'$ are
surjective. Moreover, let $Y$ has no direct summand of the form
$(P,0,0)$, where $P$ is a~projective $KI^-$-module. By \cite[Lemma
3.3]{ps}, $\overline{\psi}:Y'\to \Theta(Y)'=\Im \overline{\psi}$
is the projective cover of $\Theta(Y)'$ in $\mod (KI^-)$. Since
$\overline{\psi}f'$ is surjective and $\overline{\psi}$ is the
projective cover, the morphism $\overline{\psi}$ is essential, and
therefore $f'$ is surjective and we are done. \epv

Let $|X|$ denotes the cardinality of a~finite set $X$. Moreover,
given $KI$-modules $X$, $Y$, let ${\rm Epi}_{KI}(X,Y)$ be the set
of all surjective $KI$-homomorphisms $f:X\to Y$ and ${\rm
Ker}\,\Theta(X,Y)$ be the set of all homomorphisms $f:X\to Y$
which are in $\Ker \,\Theta$ (in the~case $X$, $Y$ are
prinjective).\vsp

{\sc Corollary 2.4.} {\it Let $K$ be a~finite field and
$X=(X',X'',\overline{\varphi})$, $Y=(Y',Y'',\overline{\psi})$ be
modules in $\prin(KI)$. If $Y$ has no direct summand of the form
$(P,0,0)$, then

$$|{\rm Epi}_{KI}(X,Y)|= |{\rm Epi}_{KI}(\Theta(X),\Theta(Y)) | \cdot
|\Ker\,\Theta (X,Y)|.$$} \vsp

{\bf Proof.} By Lemma 2.2(b) and Lemma 2.3 the functor $\Theta$
induces the surjective $K$-linear map
$$\Theta :{\rm Epi}_{KI} (X,Y) \to {\rm Epi}_{KI}
(\Theta(X),\Theta(Y))$$ by attaching to any surjective
homomorphism $f:X\to Y$ the surjective homomorphism $\Theta (f)$.
Lemma 2.3(a) finishes the proof. \epv

{\sc Lemma 2.5.} {\it Let $K$ be a~finite field and
$X=(X',X'',\overline{\varphi})$, $Y=(Y',Y'',\overline{\psi})$,
$Z=(Z',0,0)$ be modules in $\prin(KI)$. Assume that $Y$ has no
direct summand of the form $(P,0,0)$, where $P$ is a~projective
$KI^-$-module.

 \ya If there exists a~surjective homomorphism $f:X\to Y$, then there exists the~unique (up to isomorphism)
projective $KI^-$-module $U'$ such that $${\bf dim}\, U'={\bf
dim}\,X'-{\bf dim}\, Y'.$$

\yb If there is no surjective homomorphism $f:X\to Y$, then there
is no surjective homomorphism $g:X\to Y\oplus Z$.

\yc Let $U'$ be the module defined in \ya if there is a~surjective
homomorphism $f:X\to Y$ and $U'=0$ otherwise. Then $$|{\rm
Epi}_{KI}(X,Y\oplus Z)|=|{\rm Epi}_{KI}(X, Y) |\cdot |{\rm
Epi}_{KI^-}(U', Z') |\cdot | \Hom_{KI^-}(Y',Z')|.$$} \vsp

{\bf Proof.} (a) Let $f:X\to Y$ be a~surjective homomorphism and
consider $U=\Ker\, f=(U',U'',\phi)$. Since the $KI^-$ modules
$X'$, $Y'$ are projective, the $KI^-$-module $U'$ is projective.
Moreover ${\bf dim}\, U'={\bf dim}\,X'-{\bf dim}\, Y'$ and $U'$ is
uniquely determined by its dimension vector (see \cite[pp
77]{ri}).

The statement (b) is clear.

(c) If there is no surjective homomorphism $g:X\to Y$, then by (b)
the formula given in (c) is clear.

Let $g=\left[\matrix{g_1\cr g_2}\right] :X\to Y\oplus Z$ be
a~surjective homomorphism such that $g_1:X\to Y$, $g_2:X\to Z$ and
let $U=\Ker\, g_1$. It follows that $g_1,g_2$ are surjective. Note
that $X'$ may be identified with $U'\oplus Y'$, because $X'$, $Y'$
are projective $KI^-$-modules and $g_1':X'\to Y'$ is surjective
with kernel isomorphic to $U'$. Therefore the condition ${\bf
dim}\, U'={\bf dim}\,X'-{\bf dim}\, Y'$ is satisfied. By
\cite[Lemma 2.3]{kos03} there is an isomorphism of $K$-vector
spaces $\Hom_{KI}(V,Z)\simeq \Hom_{KI^-}(V',Z')$ for any
$KI$-module $V$. This isomorphism is given by $(f',f'')\mapsto f'$
and is based on the observation that $f''=0$ if $Z=(Z',0,0)$.
Therefore $g_2$ may be identified with $g_2=[g_{21},\,
g_{22}]:U'\oplus Y'\to Z'$, where $g_{21}:U'\to Z'$, $g_{22}:Y'\to
Z'$. Consider the following commutative diagram with exact rows
$$\begin{array}{ccccccccc} 0& \to & U& \hookrightarrow
& X & \longarr{g_1}{15} &Y &\to & 0 \vspace{1ex}
\\&& \mapdown{g_{21}} && \mapdown{\scriptsize{\left[\matrix{g_1\cr
g_2}\right]}} && \mapdown{\scriptsize{\id}} && \\ 0& \to & Z&
\longarr{\scriptsize{\left[\matrix{0\cr 1}\right]}}{15}  & Y\oplus
Z & \longarr{[1\,0]}{15} &Y &\to & 0\, .\end{array} $$ Since $g$
is surjective, by the Snake Lemma $g_{21}$ is surjective. So, with
any surjective $KI$-homomorphism $g:X\to Y\oplus Z$ we associate
two surjective $KI$-homomorphisms $g_1:X\to Y$, $g_{21}:U\to Z$
(identified with the surjective $KI^-$-homomorphism $g_{21}:U'\to
Z'$) and a~$KI^-$-homomorphism $g_{22}:Y'\to Z'$.

Conversely, let $g_1:X\to Y$ be a~surjective $KI$-homomorphism and
$U=\Ker\,g_1$. Note that $X'\simeq U'\oplus Y'$, because $U'$,
$X'$ and $Y'$ are projective $KI^-$-modules. Let $g_{21}:U'\to Z'$
be a~surjective $KI^-$-homomorphism and $g_{22}:Y'\to Z'$ any
$KI^-$-homomorphism. Then $g_2=[g_{21},g_{22}]:X\to Z$ is
surjective (identified with $g_2:U'\oplus Y'\to Z'$). Finally we
get a~surjective $KI$-homomorphism $g=\left[ \matrix{g_1\cr g_2}
\right]:X\to Y\oplus Z$. Indeed, let $(y,z)\in Y\oplus Z$. Let us
fix the decomposition of $X\simeq U'\oplus Y'\oplus X''$ as
a~$K$-linear space. Since $g_1$ is surjective and $g_1(U)=0$,
there exists $x_1=(0,x_1',x_1'')$ such that $g_1(x_1)=y$. Moreover
$g_{21}$ is surjective, then there exists $x_2\in U'\subseteq X$
such that $g_{21}(x_2)=z-g_{22}(x_1')$. Let $x=(x_2,x_1',x_1'')$,
therefore $g(x)=(g_1(x_1),z-g_{22}(x_1')+g_{22}(x_1'))=(y,z)$ and
lemma follows.   \epv

{\sc Lemma 2.6.} {\it Let $I$ be an~arbitrary finite poset, and
$KI$ - its incidence $K$-algebra. Let $P=\bigoplus_{i\in
I}P(i)^{n_i}$, $n_i\geq 0$, $Q=\bigoplus_{i\in I}P(i)^{m_i}$,
$m_i\geq 0$ be projective $KI$-modules. Then
$\dim_K\Hom_{KI}(P,Q)=\sum_{i\in I}(\sum_{j\preceq i}n_im_j)$. In
particular $\dim_K\Hom_{KI}(P,Q)$ is independent on the base field
$K$.}\vsp

{\bf Proof.} Let us recall that
$\dim_K\Hom_{KI}(P(i),X)=\dim_KX_i$ (see \cite[pp 68]{ri}).
Moreover $P(i)_j\simeq K$ if $i\preceq j$ in $I$ and $P(i)_j=0$
otherwise. Therefore lemma follows easily. \epv

\section{Hall polynomials for posets of finite prinjective type}
Let $I$ be a~poset of finite prinjective type and let $KI$ be its
incidence $K$-algebra. In this section we prove the~existence of
Hall polynomials for prinjective $KI$-modules. Given finite
dimensional $KI$-modules $X$, $Y$, $Z$ we define $F^Y_{Z,X}$ to be
the~number of modules $U\subseteq Y$ such that $U\simeq X$ and
$Y/U\simeq Z$.

It follows from \cite{s93}, \cite{hosim} that the Auslander-Reiten
quiver $\Gamma_I=\Gamma(\prin(KI))$ (resp. $\Gamma_{I-{\rm
sp}}=\Gamma(\modsp(KI))$) of the category $\prin(KI)$ (resp.
$\modsp(KI)$) is directed and coincides with its preprojective
component. Moreover $\Gamma_I$ and $\Gamma_{I-{\rm sp}}$ do not
depend on the base field $K$ (see \cite[Chapter 11]{s92}). Let us
recall that, by the definition, the vertices of Auslander-Reiten
quiver corresponds bijectively to the isomorphism classes of
indecomposable modules. For a~given vertex $x\in (\Gamma_I)_0$
(resp. $x\in (\Gamma_{I-{\rm sp}})_0$) we denote by $M(K,x)$
(resp. $M_{\rm sp}(K,x)$) the corresponding indecomposable
prinjective (resp. socle projective) $KI$-module. Moreover for any
function $a:(\Gamma_I)_0\to \mathbb{N}$ (resp. $a:(\Gamma_{I-{\rm
sp}})_0\to \mathbb{N}$) let $M(K,a)=\bigoplus_{x\in
(\Gamma_I)_0}M(K,x)^{a(x)}$ (resp. $M_{\rm
sp}(K,a)=\bigoplus_{x\in (\Gamma_{I-{\rm sp}})_0}M_{\rm
sp}(K,x)^{a(x)}$) (see \cite{ringel1} for details). Moreover given
a~function $a\in \CB$ we denote by $\Theta(a)\in \CB_{sp}$ the
function corresponding to the socle projective $KI$-module
$\Theta(M(a))$. It follows from \cite{s92}, \cite{s93},
\cite{hosim} and \cite{ps} that the dimension vectors $\bdim
M(K,a)$ and $\bdim M_{\rm sp}(K,a)$ depend only on the
Auslander-Reiten quiver, so they do not depend on $K$. For the
sake of simplicity we write $M(a)$ (resp. $M_{\rm sp}(a)$) instead
of $M(K,a)$ (resp. $M_{\rm sp}(K,a)$) if the base field $K$ is
known from the context. Denote by $\CB$ (resp. $\CB_{\rm sp}$) the
set of all functions $a:(\Gamma_I)_0\to \mathbb{N}$ (resp.
$a:(\Gamma_{I-{\rm sp}})_0\to \mathbb{N}$). It is clear that $\CB$
(resp. $\CB_{\rm sp}$) can be identified with the set of all
finite dimensional prinjective (resp. socle projective)
$KI$-modules. Given an~arbitrary $KI$-module $M$ we denote by
$\CS(M)$ the set of all $KI$-modules $N$ such that $\bdim N< \bdim
M$ (i.e. $\bdim N\neq \bdim M$ and $(\dim N)(i)\leq (\bdim M)(i)$
for all $i\in I$). \vsp

{\sc Lemma 3.1.} {\it Let $I$ be a~poset of finite prinjective
type. For any $a,b\in \CB$ (resp. $\overline{a},\overline{b}\in
\CB_{\rm sp}$) the natural number
$h(a,b)=\dim_K\Hom_{KI}(M(a),M(b))$ (resp.
$h(\overline{a},\overline{b})=\dim_K\Hom_{KI}(M_{\rm
sp}(\overline{a}),\allowbreak M_{\rm sp}(\overline{b}))$) does not
depend on the field $K$.} \vsp

{\bf Proof.} Since the Auslander-Reiten quivers $\Gamma_I$ and $\Gamma_{I-{\rm sp}}$ are directed,
the arguments given in \cite{ringel1} prove our lemma.\epv

For $a,b\in \CB$ (resp. $\overline{a},\overline{b}\in\CB_{sp}$) we
define polynomial $\gamma_{ab}=T^{h(a,b)}\in\mathbb{Z}[T]$ (resp.
$\gamma_{\overline{a}\overline{b}}=T^{h(\overline{a},\overline{b})}\in\mathbb{Z}[T]$).
Note that $\gamma_{ab}(|K|)=|\Hom_{KI}(M(a),M(b))|$ (resp.
$\gamma_{\overline{a}\overline{b}}(|K|)=|\Hom_{KI}(M_{sp}(\overline{a}),M_{sp}(\overline{b}))|$).
\vsp

{\sc Lemma 3.2.} {\it Let $a,b\in \CB$ and let
$\overline{a},\overline{b}\in \CB_{\rm sp}$ be such that
$\Theta(M(a))=M(\overline{a})$ and $\Theta(M(b))=M(\overline{b})$.

\ya
$|\Ker\,\Theta(M(a),M(b))|=|K|^{h(a,b)-h(\overline{a},\overline{b})}$.

\yb There exists a~polynomial $\omega _{ab}\in \mathbb{Z}[T]$ such
that for any finite field $K$ we have $\omega_{ab}(|K|)=|\Ker\,
\Theta(M(a),M(b))|$.  } \vsp

{\bf Proof.} (a) By Lemma 3.1 the natural numbers $h(a,b)$ and
$h(\overline{a},\overline{b})$ are independent on the base field
$K$. So let us fix a~finite field $K$. By Lemma 2.2(b) we have
$$|\Hom_{KI}(M(a),M(b))|=|\Hom_{KI}(M(\overline{a}),M(\overline{b}))|\cdot |\Ker\,\Theta(M(a),M(b))|.$$
To finish the prove of (a) we have only to observe that
$|\Hom_{KI}(M(a),M(b))|=|K|^{h(a,b)}$ and
$|\Hom_{KI}(M(\overline{a}),M(\overline{b}))|=|K|^{h(\overline{a},\overline{b})}$.

(b) Put $\omega_{ab}=T^{h(a,b)-h(\overline{a},\overline{b})}$.
Then (b) follows from (a). \epv

{\sc Theorem 3.3.} {\it Let $I$ be a~poset of finite prinjective
type and let $a\in \CB$ (resp. $\overline{a}\in\CB_{\rm sp}$).
There exists a~monic polynomial $\alpha_{a}\in \mathbb{Z}[T]$
(resp. $\alpha_{\overline{a}}\in \mathbb{Z}[T]$) such that for any
finite field $K$
$$|\Aut_{KI}(M(a))|=\alpha_{a}(|K|), \;\; (\mbox{resp. } |\Aut_{KI}(M_{\rm sp}(\overline{a}))|=\alpha_{\overline{a}}(|K|)). $$ }

{\bf Proof.} We may follow the proof given in \cite{ringel1}.
This theorem  also follows from \cite[Proposition 2.1]{peng}.\epv

Given functions $x,y,z\in \CB\cup \CB_{sp}$, etc., for the sake of
simplicity, we denote by capital letters $X$, $Y$, $Z$, etc. the
$KI$-modules $M(K,x)$, $M_{sp}(K,x)$, $M(K,y)$, $M(K,z)$,
respectively. However we should remember that $KI$-modules are
identified with functions from the sets $\CB$, $\CB_{sp}$ and
depend on the base field $K$. Moreover given a~function $x\in \CB$
we denote by $\Theta(x)$ the function in $\CB_{sp}$ corresponding
to the module $\Theta(X)$.\vsp

{\sc Lemma 3.4.} {\it Let $I$ be a~poset of finite prinjective
type. Let $x,y\in\CB_{sp}$. There exist polynomials $\sigma
^y_x,\eta^y_x,\mu^y_x,\varepsilon^y_x\in \mathbb{Z}[T]$ such that
for any finite field $K$:

$\sigma ^y_x(|K|)$ equals the number of submodules $U\subseteq Y$,
such that $U\simeq X$,

$\eta ^y_x(|K|)$ equals the number of submodules $U\subseteq Y$,
such that $Y/U\simeq X$,

$\mu^y_x(|K|)$ equals the number of injective homomorphisms $X\to
Y$,

$\varepsilon^y_x(|K|)$ equals the number of surjective
homomorphisms $Y\to X$.} \vsp

{\bf Proof.} One can prove this lemma by
developing Ringel's arguments given in \cite{ringel1}. For the
convenience of the reader we outline the proof.

If ${\bf dim}\, X\nleqslant {\bf dim}\, Y$, we set
$\sigma^y_x=0=\eta^y_x$.

Let ${\bf dim}\, X\leqslant {\bf dim}\, Y$. We apply induction on
${\bf dim}\, Y$. If ${\bf dim}\, Y=0$, then $X=0=Y$ and
$\sigma^y_x=1=\eta^y_x$. Let $Y\neq 0$ and we start with induction
on ${\bf dim}\, X$. Define two polynomials
$\mu^y_x=\gamma_{xy}-\sum_{U\in
\CS(X)}\eta^x_u\alpha_u\sigma^y_u$,
$\;\varepsilon^y_x=\gamma_{yx}-\sum_{U\in
\CS(X)}\eta^y_u\alpha_u\sigma^x_u$.
Since the category $\modsp(KI)$ is closed under submodules,
we may assume that $U$ arising in these sums is
socle projective, because otherwise $\sigma^y_u=0=\sigma^x_u$.
Moreover these sums are finite, because the poset $I$ is of finite prinjective type. All summands on the right side are defined by induction
hypothesis.

We claim that $\eta^x_u\alpha_u\sigma^y_u(|K|)$ equals the number
of morphisms $f:X\to Y$ such that $\Im f\simeq U$. Indeed, for
a~given submodule $V\subseteq X$ such that $X/V\simeq U$ we fix
a~surjective homomorphism $g_V:X\to U$ with $\Ker\, g_V=V$.
Similarly, if $W\subseteq Y$ is a~submodule such that $W\simeq U$,
we fix an~injective homomorphism $h_W:U\to Y$ with $\Im h_W=W$.
Homomorphisms $X\to Y$ with kernel $V$ and image $W$ correspond
bijectively to automorphisms of $U$. This bijection is given by
attaching to any automorphism $f:U\to U$ the following
homomorphism $X\to Y$:
$$X\longarr{g_V}{30}U\longarr{f}{30}U\longarr{h_W}{30}Y.$$
A~homomorphism $X\to Y$ is injective if and only if its image is
not isomorphic to any $U$ with ${\bf dim}\,  U< {\bf dim}\, X$.
Therefore $\mu^y_x(|K|)$ is the number of injective homomorphisms
$X\to Y$. Dually, $\varepsilon^y_x(|K|)$ is the number of
surjective homomorphisms $Y\to X$.

Note that for all finite fields $K$, $\mu^y_x(|K|)(\alpha_x(|K|))^{-1}$ equals the number of submodules
$U\subseteq Y$ with $U\simeq X$ and therefore it is an~integer. By
\cite[page 441]{ringel1} the polynomial $\alpha_x$ divides $\mu^y_x$ in $\mathbb{Z}[T]$.
Similarly, $\alpha_x$ divides $\varepsilon^y_x$ in $\mathbb{Z}[T]$. We put
$\sigma^y_x=\mu^y_x(\alpha_x)^{-1}$ and
$\eta^y_x=\varepsilon^y_x(\alpha_x)^{-1}$. This finishes the
proof.\epv

{\sc Lemma 3.5.} {\it Let $I$ be an~arbitrary poset and let $X$,
$Y$ be projective $KI$-modules there exist polynomials $\eta^y_x,
\varepsilon^y_x\in \mathbb{Z}[T]$ such that for any finite field
$K$:

$\eta ^y_x(|K|)$ equals the number of submodules $U\subseteq Y$
such that $Y/U\simeq X,$

$\varepsilon ^y_x(|K|)$ equals the number of surjective
homomorphisms $Y\to X$.} \vsp

{\bf Proof.} Let $X$, $Y$, $Z$ be $KI$-modules. By \cite[Section 4]{riedtmann},
the number of submodules $U\subseteq Y$, such that
$U\simeq Z$ and $Y/U\simeq X$, equals
$$F^Y_{X,Z}=\frac{|\Ext^1_{KI}(X,Z)_Y||\Aut_{KI}(Y)|}{|\Aut_{KI}(Z)||\Aut_{KI}(X)|
|\Hom_{KI}(Z,X)|},\leqno(*)$$ where $\Ext^1_{KI}(X,Z)_Y$ is the
set of all exact sequences in $\Ext^1_{KI}(X,Z)$ with the middle
term $Y$. Let us assume that $Y$ and $X$ are projective
$KI$-modules. Let us fix a~submodule $Z\subseteq Y$ such that
$Y/Z\simeq X$. Since the category of projective modules is closed
under kernels of surjective homomorphisms, the submodules
$U\subseteq Y$ with $Y/U\simeq X$ are projective. Moreover
$U\simeq Z$, because any exact sequence $0\to U\to Y\to X\to 0$
splits. Therefore $F^Y_{X,Z}$ equals the number of submodules
$U\subseteq Y$ such that $Y/U\simeq X$. Note also that
$\Ext^1_{KI}(X,Z)=0$ and therefore $|\Ext^1_{KI}(X,Z)_Y|=1$. By
Lemma 2.6 the number $h(z,x)=\dim_K\Hom_{KI}(Z,X)$ is independent
on the base field $K$ and the number of $KI$-homomorphisms $f:Z\to
X$ equals $\gamma_{z,x}(|K|)$. We define
$$\eta^y_x=\frac{\alpha_y}{\alpha_z\alpha_x\gamma_{z,x}}.$$ By
Theorem 3.3 and $(*)$, $F^Y_{X,Z}=\eta^y_x(|K|)$ for any finite
field $K$. Then the number
$$\alpha_z(|K|)\alpha_x(|K|)\gamma_{z,x}(|K|)$$ divides
$\alpha_y(|K|)$ for infinitely many finite fields $K$. Since the
polynomial $\alpha_z\alpha_x\gamma_{z,x}$ is monic, it follows
from \cite[page 441]{ringel1} that it divides the polynomial
$\alpha_y$ in $\mathbb{Z}[T]$ and therefore $\eta^y_x\in
\mathbb{Z}[T]$. Consequently $\eta^y_x (|K|)$ equals the number of
submodules $U\subseteq Y$ such that $Y/U\simeq X$.

We put $\varepsilon^y_x=\eta^y_x\alpha_x\in \mathbb{Z}[T]$. Note
that $\varepsilon^y_x(|K|)$ equals the number of surjective
homomorphisms $f:Y\to X$. This finishes the proof.\epv

{\sc Corollary 3.6.} {\it  Assume that $I$ is of finite
prinjective type and $x,y\in \CB$. There exists a~polynomial
$\varepsilon^y_x\in \mathbb{Z}[T]$ such that for any field $K$:

$\varepsilon ^y_x(|K|)={\rm Epi}_{KI}(Y, X)$.} \vsp

{\bf Proof.} If there is no surjective homomorphism $f:Y\to X$ for
any field $K$, we put $\varepsilon^y_x=0$. Otherwise, let
$X=\overline{X}\oplus Z$, where $Z=(P,0,0)$ with projective
$KI^-$-module $P$ and $\overline{X}$ has no direct summand of the
form $(P,0,0)$. Then $\Theta(X)=\Theta(\overline{X})$. In our case
there exists a~surjective homomorphism $f:Y\to X$ for some field
$K$. Let $U'\simeq Y'/\overline{X}'$ be the unique (up to
isomorphism) projective $KI^-$-module such that $\bdim\,
U'=\bdim\, Y'-\bdim\, \overline{X}'$. By Lemma 3.4, there exists
a~polynomial $\varepsilon^{\Theta(y)}_{\Theta(x)}\in
\mathbb{Z}[T]$ such that $\varepsilon
^{\Theta(y)}_{\Theta(x)}(|K|)$ equals the number of surjective
homomorphisms $\Theta(Y)\to \Theta(X)$. By Lemma 3.5, there exists
a~polynomial $\varepsilon ^{u'}_{z'}\in \mathbb{Z}[T]$ such that
$\varepsilon ^{u'}_{z'}(|K|)$ equals the number of surjective
homomorphisms $U'\to Z'$. Put
$$\varepsilon^y_x= \varepsilon ^{\Theta(y)}_{\Theta(x)}\cdot T^{h(y,\overline{x})-h(\Theta(y),\Theta(x))}
\cdot T^{h(x',z')}\cdot \varepsilon^{u'}_{z'}.$$  By Corollary
2.4, Lemma 2.5 and Lemma 3.2, $\varepsilon^y_x$ is the required
polynomial.\epv

{\sc Corollary 3.7.} {\it  Let $I$ be a~poset of finite
prinjective type
 and let $x,y\in \CB$.
There exists a~polynomial $\eta^y_x\in \mathbb{Z}[T]$ such that
for any finite field $K$:

$\eta ^y_x(|K|)$ equals the number of submodules $U\subseteq Y$,
such that $Y/U\simeq X$.}\vsp

{\bf Proof.} By Corollary 3.6, there exists a~polynomial
$\varepsilon^y_x\in \mathbb{Z}[T]$ such that $\varepsilon
^y_x(|K|)={\rm Epi}_{KI}(Y,X)$ for any finite field $K$. Note
that, for any finite field $K$, the number
$\varepsilon^y_x(|K|)\cdot \alpha_x(|K|)^{-1}$ is an~integer,
because it counts the number of submodules $U\subseteq Y$ such
that $Y/U\simeq X$. Since $\alpha_X$ is a~monic polynomial, it
follows from \cite[page 441]{ringel1} that $\alpha_x$ divides
$\varepsilon^y_x$ in $\mathbb{Z}[T]$. Therefore
$\eta^y_x=\varepsilon^y_x\cdot \alpha_x^{-1}\in \mathbb{Z}[T]$ is
the required polynomial.\epv

{\sc Theorem 3.8.} {\it  Let $I$ be a~poset of finite prinjective
type and $x$, $y$, $z$ be functions in $\CB$ (resp.
$\overline{x},\overline{y},\overline{z}\in\CB_{sp}$). There exist
polynomials $\varphi ^y_{xz}\in \mathbb{Z}[T]$ (resp. $\varphi
^{\overline{y}}_{\overline{x}\overline{z}}\in \mathbb{Z}[T]$) such
that for any finite field $K$:

$\varphi ^y_{xz}(|K|)=F^Y_{XZ}$ (resp. $\varphi
^{\overline{y}}_{\overline{x}\overline{z}}(|K|)=F^{\overline{Y}}_{\overline{X}\overline{Z}}$).
}\vsp

{\bf Proof.} We prove this theorem developing arguments given in
\cite{ringel1} and facts proved in Sections 2 and 3.

If ${\bf dim}\, Y\neq {\bf dim}\, Z+{\bf dim}\, X$ we put $\varphi
^y_{xz}=0$. Let ${\bf dim}\, Y= {\bf dim }\, Z+{\bf dim}\, X$. We
apply induction on ${\bf dim}\, Z$. If ${\bf dim}\, Z=0$
we put $\varphi^x_{x0}=1$ and $\varphi^y_{x0}=0$ if
$X\not\simeq Y$.

Assume that $Z\neq 0$ and $Z=U_1\oplus U_2$, where $U_1\neq 0$,
$U_1\simeq W^m$, $W$ is indecomposable, $W$ is not a~direct
summand of $U_2$ and no indecomposable direct summand of $U_2$ is
a~predecessor of $W$ in $\Gamma_I$ (resp. $\Gamma_{I-\rm sp}$ in
the~"socle projective" case).

Let us consider two cases:

{\bf Case 1.} $U_2\neq 0$. We define
$$\varphi^y_{xz}=\sum_d\varphi ^d_{xu_1}\varphi^y_{du_2},
$$ where the sum runs over all modules $D$ such that
${\bf dim}\, D={\bf dim}\, X+{\bf dim}U_1$. Note that this sum is
finite and runs over prinjective modules (resp. socle projective
modules), because the category of prinjective modules (resp. socle
projective modules) is closed under extensions and the poset $I$
is of finite prinjective type. Moreover the right
side is already defined by induction hypothesis. One can
prove that $\varphi^y_{xz}(|K|)=F^Y_{XZ}$ (see \cite{ringel1}).

{\bf Case 2.} $U_2=0$. We define
$$\varphi^y_{xz}=\eta^y_x-\sum_{d\not\simeq z}\varphi^y_{xd},
$$ where $d$ runs over all modules such that
${\bf dim}\, D={\bf dim}\, Z$. Since the category of prinjective
is closed under kernels of epimorphisms and the category of socle
projective modules is closed under submodules, we may assume that
the modules $D$ are prinjective (resp. have projective socle).
Note that $D$ is not a~direct power of indecomposable, because $Z$
is a~direct power of indecomposable, $Z\not\simeq D$ and $\bdim\,
Z=\bdim\, D$ (see \cite[IX.2.1]{ars}). Therefore the polynomials
$\varphi^y_{xd}$ are defined in Case 1. The polynomials $\eta^y_x$
are defined in Corollary 3.7 for prinjective modules and in Lemma
3.4 for socle projective modules. It is clear that $\varphi
^y_{xz}(|K|)=F^Y_{XZ}$  and this finishes the proof.\epv

The polynomials $\varphi^y_{xz}$ are called {\bf Hall
polynomials}.\vsp

%
%
%
%
%
%
%
%

In the last chapter we present consequences of the~existence of
Hall polynomials for prinjective modules.\vsp

\section{Prinjective Ringel-Hall algebras} We denote by
$\CH_{prin}(I)$ the free $\mathbb{Q}(T)$-module with basis
$\{u_x\}_{x\in \CB}$, indexed by the elements of the set $\CB$.
$\CH_{prin}(I)$ is equipped with a~multiplication defined by the
formula:
$$u_{x_1}u_{x_2}=\sum_{x\in \CB}\varphi^x_{x_1x_2}u_x. $$  Note
that this sum is finite, because the poset $I$ is of finite
prinjective type and $\varphi^x_{x_1,x_2}\neq 0$ only
 if $\bdim X=\bdim X_1+\bdim X_2$. By
\cite[Proposition 4]{ringel1}, $\CH_{prin}(I)$ is an~associative
ring and the element $u_0$ is the identity element of
$\CH_{prin}(I)$. By the results of Section 3 this ring depends
only on the poset $I$. We call $\CH_{prin}(I)$ the {\bf
prinjective generic Ringel-Hall algebra} for the poset $I$. \vsp

{\bf Concluding remarks.} (1) In the forthcoming paper
\cite{kos05b} description of $\CH_{prin}(I)$ by generators and
relations is given. Moreover in \cite{kos05b} we show connections
of the prinjective Ringel-Hall algebra with Lie algebras and
Kac-Moody algebras.

(2) In \cite{kos04} the existence of generic extensions for
prinjective modules over posets of finite prinjective type is
proved. It would be interesting to find connections between the
monoid of generic extensions of prinjective modules and some
specialization of prinjective Ringel-Hall algebra. Such
a~connection, for Dynkin quivers, one ca find in
\cite{reineke2001}.

(3) In the paper \cite{kos05b} generators of prinjective
Ringel-Hall algebra are given. Most of these generators are in the
kernel of the functor $\Theta$. We can't see natural candidates
for generators in the "socle projective case", therefore the
category of prinjective modules is more convenient in our
considerations.

(4) In \cite{kub} the existence of Hall polynomials for
representations of finite type bisected posets is proved. However,
in our case, it solves only the problem of the existence of Hall
polynomials for socle projective modules over posets of finite
prinjective type with exactly one maximal element.

%

\end{document}